# Explicit expressions for the variogram of first–order intrinsic autoregressions

**Tibor K. Pogány**

*Faculty of Maritime Studies, University of Rijeka, 51000 Rijeka, Studentska 2, Croatia*
*e-mail:* `poganj@pfri.hr`
*e-mail:* `www.pfri.hr/ poganj`

**Saralees Nadarajah**

*School of Mathematics, University of Manchester, Manchester M13 9PL, U.K.*
*e-mail:* `mssbbbn2@manchester.ac.uk`

**Abstract:** Exact and explicit expressions for the variogram of first–order intrinsic autoregressions have not been known. Various asymptotic expansions and approximations have been used to compute the variogram. In this note, an exact and explicit expression applicable for all parameter values is derived. The expression involves Appell's hypergeometric function of the fourth kind. Various particular cases of the expression are also derived.

## 1. Introduction

Let $\{X_{u,v} : u, v \in \mathbb{Z}\}$ be a homogeneous first–order intrinsic autoregression on the two–dimensional rectangular lattice $\mathbb{Z}^2$ [11], [3] with generalized spectral density function

$$f(w, \eta) = \kappa \left(1 - 2a \cos w - 2b \cos \eta\right)^{-1}$$

for $w \in (-\pi, \pi)$, $\eta \in (-\pi, \pi)$ and the conditional expectation structure

$$\mathsf{E}\left(X_{u,v} \mid \cdots\right) = a\left(x_{u-1,v} + x_{u+1,v}\right) + b\left(x_{u,v-1} + x_{u,v+1}\right) \qquad (1)$$

where $a > 0$, $b > 0$, $a + b = \frac{1}{2}$ and

$$\mathrm{var}\left(X_{u,v} \mid \cdots\right) = \kappa \,.$$

We can assume without loss of generality that $\kappa = 1$ and that $\{X_{u,v}\}$ is Gaussian, [4]. It follows that the difference $X_{u,v} - X_{u+s,v+t}$ has a well defined distribution with zero mean and lag $(s,t)$ variogram

$$\nu_{st}(a,b) := \frac{1}{\pi^2} \int_0^\pi \int_0^\pi \frac{1 - \cos(sx)\cos(ty)}{1 - 2a \cos x - 2b \cos y} \, \mathrm{d}x \mathrm{d}y \,. \qquad (2)$$

The computation of (2), both analytically and numerically, has been a subject of considerable interest. The symmetric case $a = b = 1/4$ was considered by McCrea and Whipple [13], Spitzer [18] and Besag and Kooperberg [3]. Besag and Mondal [4] derived explicit expressions for $\nu_{s,s}$ and $\nu_{s,0}$ for the asymmetric case with the latter expressed as a finite sum of incomplete beta functions. Duffin and Shaffer [9, Theorem 4] and Besag and Mondal [4, Theorem 2] provided asymptotic expansions for $\nu_{s,t}$ in terms of $r = \sqrt{4bs^2 + 4at^2}$ for the symmetric and asymmetric cases. Their approach was to exploit the recurrence equation with respect to integer time variables $s$ and $t$ [4, Eq. (4)]. See also [2] and [11].

The aim of this note is to derive an explicit expression for (2) for general $s, t, a$ and $b$. Our approach is quite different — we attack the problem varying $a$ and $b$ to find the explicit formula for $\nu_{st}(a, b)$.

The expression given in Section 2 involves Appell's hypergeometric function of the fourth kind [10, page 1008] defined by

$$F_4\big[\alpha, \beta; \gamma, \gamma'; x, y\big] = \sum_{j=0}^{\infty} \sum_{k=0}^{\infty} \frac{(\alpha)_{j+k}(\beta)_{j+k}}{(\gamma)_j\,(\gamma')_k} \frac{x^j}{j!} \frac{x^k}{k!} \tag{3}$$

for $\sqrt{|x|} + \sqrt{|y|} < 1$, where $(w)_\ell := w(w+1)\cdots(w+\ell-1)$ denotes Pochhammer symbol with $(w)_0 \equiv 1$. Various particular cases of the general expression, involving simpler functions, are derived in Section 3.

A transformation formula between $F_4$ and the Appell's hypergeometric function $F_2$ defined by

$$F_2\big[\alpha, \beta, \beta'; \gamma, \gamma'; x, y\big] := \sum_{j=0}^{\infty} \sum_{k=0}^{\infty} \frac{(\alpha)_{j+k}(\beta)_j(\beta')_k}{(\gamma)_j\,(\gamma')_k} \frac{x^j}{j!} \frac{x^k}{k!} \qquad |x| + |y| < 1$$

is [17, §5.4]

$$F_4\big[\alpha, (\alpha+1)/2; \gamma + \tfrac{1}{2}, \gamma' + \tfrac{1}{2}; x^2, y^2\big]$$
$$= (1 + x + y)^{-\alpha} F_2\Big[\alpha, \gamma, \gamma'; 2\gamma, 2\gamma'; \frac{2x}{x+y+1}, \frac{2y}{x+y+1}\Big].$$

This formula has been proved earlier (in equivalent form) in [1, p. 11, Eq. (3.1)], see [21, Eq. (175), §9.4.] too. In–built numerical routines for the computation of $F_2$ are available, see [8, §3.1.3] and especially [7, §2]. Reduction procedures for $F_4$ to lower order analytical expressions have been developed by Niukkanen [14], Paramonova and Niukkanen [15] and references therein.

## 2. Main result

Theorem 1 expresses (2) in terms of Appell's hypergeometric function of the fourth kind defined in (3). It applies for $|a| + |b| < \tfrac{1}{2}$. The case $|a| + |b| = \tfrac{1}{2}$ is considered by Theorem 2.

**Theorem 1.** *For all $s, t \in \mathbb{N}_0$, $|a| + |b| < \frac{1}{2}$ we have*

$$\nu_{st}(a,b) = F_4\left[\tfrac{1}{2}, 1; 1, 1; 4a^2, 4b^2\right] - \binom{s+t}{s} a^s b^t$$
$$\times F_4\left[\tfrac{1}{2}(s+t+1), \tfrac{1}{2}(s+t)+1; s+1, t+1; 4a^2, 4b^2\right]. \quad (4)$$

*Proof.* By equation **3.915**(2) in [10], we can write

$$\nu_{st}(a,b) = \frac{1}{\pi^2} \int_0^\pi \int_0^\pi \int_0^\infty \left(1 - \cos(sx)\cos(ty)\right) e^{-(1-2a\cos x - 2b\cos y)z} \, dxdydz$$
$$= \int_0^\infty e^{-x} \left( J_0(-2ax\mathrm{i}) J_0(-2bx\mathrm{i}) - \mathrm{i}^{s+t} J_s(-2ax\mathrm{i}) J_t(-2bx\mathrm{i}) \right) dx \quad (5)$$

where $\mathrm{i} = \sqrt{-1}$ and $J_s(x)$ denotes the Bessel function of the first kind of order $s$. Setting

$$\mathcal{I}_{st}(a,b) := \mathrm{i}^{s+t} \int_0^\infty e^{-x} J_s(-2ax\mathrm{i}) J_t(-2bx\mathrm{i}) \, dx$$

we can rewrite (5) as

$$\nu_{st}(a,b) = \mathcal{I}_{00}(a,b) - \mathcal{I}_{st}(a,b). \quad (6)$$

Consider now for instance the Laplace–transform formula **3.12.15**(20) in [16], *viz*

$$\int_0^\infty e^{-pt} t^\lambda J_\mu(At) J_\nu(Bt) \, dt = \frac{A^\mu B^\nu}{2^{\mu+\nu} p^{\lambda+\mu+\nu+1}} \frac{\Gamma(\lambda+\mu+\nu+1)}{\Gamma(\mu+1)\Gamma(\nu+1)}$$
$$\times F_4\left[\tfrac{1}{2}(\lambda+\mu+\nu+1), \tfrac{1}{2}(\lambda+\mu+\nu+2); \mu+1, \nu+1; -\frac{A^2}{p^2}, -\frac{B^2}{p^2}\right] \quad (7)$$

when $\Re\{\lambda+\mu+\nu\} > -1$, $\Re\{p\} > |\Im\{A\}| + |\Im\{B\}|$. Setting $p = 1$, $\lambda = 0$, $A = -2a\mathrm{i}$, $B = -2b\mathrm{i}$ in (7) one deduces

$$\mathcal{I}_{st}(a,b) = \binom{s+t}{s} a^s b^t F_4\left[\tfrac{1}{2}(s+t+1), \tfrac{1}{2}(s+)+1; s+1, t+1; 4a^2, 4b^2\right]. \quad (8)$$

The result of the theorem follows from (6) and (8). □

**Remark 1.** *The explicit expression in (4) has some applicability. Although one assumes that $|a| + |b| = \frac{1}{2}$ in (1) the case $|a| + |b| < \frac{1}{2}$ has been of interest. For instance, see [3, Example 3.3] quoting a result by Kempton & Howes on a four–neighbour auto–normal scheme with the maximum likelihood estimates $\widehat{a} = 0.4848$ and $\widehat{b} = 0.0132$, i.e. $\widehat{a} + \widehat{b} = 0.4980 < \frac{1}{2}$.*

**Theorem 2.** *For all $a \in \left(0, \tfrac{1}{2}\right)$ we have*

$$\nu_{st}\left(a, \tfrac{1}{2} - a\right) = \lim_{\theta \to 0+} \nu_{st}\left(a\sqrt{1-\theta}, \left(\tfrac{1}{2} - a\right)\sqrt{1-\theta}\right). \quad (9)$$

*Proof.* Consider the Laplace transform (7) for $p = (1-\theta)^{-1/2}$ and $\lambda = 0$. The resulting Appell's function $F_4$ obviously converges on the edge $a + b = \frac{1}{2}$ of the convergence region. Considering that $\nu_{st}$ is a difference of two $F_4$ terms, the result follows by Abel's summation method. □

**Remark 2.** *The limit in* (9) *gives an approximation formula, viz*

$$\nu_{st}(a, \tfrac{1}{2} - a) \approx F_4\left[\tfrac{1}{2}, 1; 1, 1; 4a^2(1-\theta), 4\left(\tfrac{1}{2} - a\right)^2(1-\theta)\right]$$
$$- \binom{s+t}{s} \frac{a^s\left(\tfrac{1}{2} - a\right)^t}{(1-\theta)^{-(s+t)/2}} F_4\left[\tfrac{1}{2}(s+t+1), \tfrac{1}{2}(s+t)+1;\right.$$
$$\left. s+1, t+1; 4a^2(1-\theta), 4\left(\tfrac{1}{2} - a\right)^2(1-\theta)\right], \quad (10)$$

*where the quality of the approximation is controlled by suitably chosen small $\theta > 0$.*

## 3. The symmetric case $a = b = 1/4$

Here, we present technical details to calculate $\nu_{s,t}\left(\tfrac{1}{4}, \tfrac{1}{4}\right)$. Using Burchnall–formula [21, §9.4, Eq. (149)], we can transform

$$F_4[\alpha, \beta; \gamma, \gamma'; x, x] = {}_4F_3\left[\begin{array}{c}\alpha, \ \beta, \ \tfrac{1}{2}(\gamma+\gamma'), \ \tfrac{1}{2}(\gamma+\gamma'-1) \\ \gamma, \ \gamma', \ \gamma+\gamma'-1\end{array} ; 4x\right].$$

The asymptotics of generalized hypergeometric series ${}_{p+1}F_p, p \geq 3$ have been studied by Bühring and Srivastava [6], Saigo and Srivastava [19] and by A.K. Srivastava [20]. Since we are faced with

$$_4F_3\left[\begin{array}{c}\tfrac{1}{2}(s+t+1), \ \tfrac{1}{2}(s+t)+1, \ \tfrac{1}{2}(s+t)+1, \ \tfrac{1}{2}(s+t+1) \\ s+1, \ t+1, \ s+t+1\end{array} ; 1\right], \quad (11)$$

we have a hypergeometric term

$$_4F_3\left[\begin{array}{c}\alpha_1, \ \alpha_2, \ \alpha_3, \ \alpha_4 \\ \beta_1, \ \beta_2, \ \beta_3\end{array} ; 1\right]$$

that is zero–balanced, i.e. when $\alpha_1 + \alpha_2 + \alpha_3 + \alpha_4 = \beta_1 + \beta_2 + \beta_3$. To account for the asymptotics of this ${}_4F_3[\cdots; 1-\theta]$ as $\theta \to 0+$, consider equation (4.2) in [6] for $p = 3$:

$$_4F_3\left[\begin{array}{c}\alpha_1, \ \alpha_2, \ \alpha_3, \ \alpha_4 \\ \beta_1, \ \beta_2, \ \beta_3\end{array} ; 1-\theta\right] = \frac{1}{\mathbf{\Gamma}}\left\{L + B - \ln\theta\right\} + O(\theta) + O(\theta\ln\theta), \quad (12)$$

valid for $\theta \to 0+$. Here

$$\mathbf{\Gamma} = \frac{\Gamma(\alpha_1)\Gamma(\alpha_2)\Gamma(\alpha_3)\Gamma(\alpha_4)}{\Gamma(\beta_1)\Gamma(\beta_2)\Gamma(\beta_3)}, \quad (13)$$

and
$$L = -2\gamma - \Psi(\alpha_1) - \Psi(\alpha_2),$$
where $\gamma = -\Psi(1)$ denotes the Euler–Mascheroni constant and $\Psi(x) = \frac{\mathrm{d}}{\mathrm{d}x}\ln\Gamma(x)$ denotes the digamma function. The $B$ in (12) can be represented by various forms, see (4.6), (4.7), and (4.12) in [6]. An example is [6, Eq. (4.7)]:

$$B = \sum_{k=1}^{\infty} \frac{(\beta_3 + \beta_1 - \alpha_4 - \alpha_3)_k(\beta_3 + \beta_2 - \alpha_4 - \alpha_3)_k}{k(\alpha_1)_k(\alpha_2)_k}$$
$$\times {}_3F_2\left[\begin{array}{c} \beta_3 - \alpha_3,\ \beta_3 - \alpha_4,\ -k \\ \beta_3 + \beta_1 - \alpha_4 - \alpha_3,\ \beta_3 + \beta_2 - \alpha_4 - \alpha_3 \end{array}; 1\right]. \qquad (14)$$

An advantage of this formula is that it is a single infinite series of hypergeometric terms and the series for each hypergeometric term is finite because of the $-k$ in the numerator. Combining (11), (13) and (14), we can write $L_{st} := L$, $B_{st} := B$ and $\mathbf{\Gamma}_{st} := \mathbf{\Gamma}$ as

$$B_{st} = \sum_{k=1}^{\infty} \frac{4^k\left(s+\frac{1}{2}\right)_k\left(t+\frac{1}{2}\right)_k}{k(s+t+1)_{2k}} {}_3F_2\left[\begin{array}{c} \frac{1}{2}(s+t),\ \frac{1}{2}(s+t+1),\ -k \\ s+\frac{1}{2},\ t+\frac{1}{2} \end{array}; 1\right] \quad (15)$$

and
$$\mathbf{\Gamma}_{st} = \mathbf{\Gamma} = \binom{s+t}{s}\frac{\pi}{4^{s+t}}$$
respectively. Note for instance that
$$B_{ss} = \Psi(s+1) - \Psi\left(s+\tfrac{1}{2}\right) \qquad s \in \mathbb{N}_0, \qquad (16)$$
see the Appendix. Finally, routine calculations show that:

$$\nu_{st}\left(\tfrac{\sqrt{1-\theta}}{4},\tfrac{\sqrt{1-\theta}}{4}\right) = \frac{\ln 4 + B_{00} - \ln\theta}{\pi}$$
$$- \binom{s+t}{s}\frac{(1-\theta)^{(s+t)/2}}{4^{s+t}\mathbf{\Gamma}_{st}}\left\{L_{st} + B_{st} - \ln\theta\right\} + O(\theta) + O(\theta\ln\theta)$$
$$= \frac{1}{\pi}\left\{\ln 16 + B_{00} + 2\gamma + \Psi\left(\tfrac{1}{2}(s+t+1)\right) + \Psi\left(\tfrac{1}{2}(s+t)+1\right) - B_{st}\right\}$$
$$+ O(\theta) + O(\theta\ln\theta) \qquad \theta \to 0_+. \qquad (17)$$

for all $(s,t) \in \mathbb{N}_0^2$.

**Theorem 3.** *For all $s,t \in \mathbb{N}_0$ we have*

$$\nu_{st}\left(\tfrac{1}{4},\tfrac{1}{4}\right) = \frac{1}{\pi}\left\{\ln 4 + 2\sum_{k=1}^{s+t}\frac{1}{k} - B_{st}\right\}. \qquad (18)$$

*Proof.* In (17) assume that $s+t$ is even. So, applying properties of digamma function $\Psi$ to expression (17), we conclude

$$\Psi\left(\tfrac{1}{2}(s+t+1)\right) + \Psi\left(\tfrac{1}{2}(s+t)+1\right) = -2\gamma - \ln 4 + \sum_{k=0}^{\frac{1}{2}(s+t)-1}\frac{1}{k+\frac{1}{2}} + \sum_{k=1}^{\frac{1}{2}(s+t)}\frac{1}{k}.$$

For odd $s+t$ repeat this procedure. □

**Corollary 3.1.** *For all $s \in \mathbb{N}_0$,*

$$\nu_{ss}\left(\tfrac{1}{4}, \tfrac{1}{4}\right) = \frac{4}{\pi} \sum_{k=0}^{s-1} \frac{1}{2k+1}. \tag{19}$$

*Proof.* Consider (18) for $s = t$. By (16) it follows

$$\nu_{s,s}\left(\tfrac{1}{4}, \tfrac{1}{4}\right) = \frac{2}{\pi}\left\{\ln 4 + \gamma + \Psi\left(s + \tfrac{1}{2}\right)\right\} = \frac{2}{\pi}\left\{\ln 4 + \gamma + \Psi\left(\tfrac{1}{2}\right) + \sum_{k=0}^{s-1} \frac{1}{k + \tfrac{1}{2}}\right\},$$

so the result. □

**Remark 3.** *The equation* (19) *is well–known in the literature, see* [3, Example 3.2]*, where the authors refer back to* [18, p. 148]*.*

### Acknowledgement

The recent investigation was supported in part by Research Project No. 112-2352818-2814 of Ministry of Sciences, Education and Sports of Croatia.

### Appendix

Here we prove (16). Using the transformation given by Bühring [5, Eq. (4.1)], one can express $B_{st}$ in (15) as

$$B_{st} = \sum_{k=1}^{\infty} \frac{4^k \left(s + \tfrac{1}{2}\right)_k \left(\tfrac{1}{2}(t - s + 1)\right)_k}{k \left(\tfrac{1}{2}(s + t + 1)\right)_k \left(\tfrac{1}{2}(s + t) + 1\right)_k} \,_3F_2\!\left[\begin{array}{c}\tfrac{1}{2}(s+t),\ \tfrac{1}{2}(s-t),\ -k \\ s + \tfrac{1}{2},\ \tfrac{1}{2}(s - t + 1) - k\end{array}; 1\right].$$

If $s = t$ then the hypergeometric term reduces to 1, so

$$\begin{aligned}
B_{ss} &= \sum_{k=1}^{\infty} \frac{\left(\tfrac{1}{2}\right)_k}{k\,(s+1)_k} = \int_0^1 \sum_{k=1}^{\infty} \frac{\left(\tfrac{1}{2}\right)_k x^{k-1}}{(s+1)_k}\,\mathrm{d}x \\
&= \frac{1}{2(s+1)} \int_0^1 \sum_{k=1}^{\infty} \frac{\left(\tfrac{3}{2}\right)_{k-1}(1)_{k-1}}{(s+2)_k} \frac{x^{k-1}}{(k-1)!}\,\mathrm{d}x \\
&= \frac{1}{2(s+1)} \int_0^1 {}_2F_1\!\left[\begin{array}{c}\tfrac{3}{2},\ 1 \\ s + 2\end{array}; x\right] \mathrm{d}x \\
&= \frac{1}{2(s+1)} \,_3F_2\!\left[\begin{array}{c}\tfrac{3}{2},\ 1,\ 1 \\ s + 2,\ 2\end{array}; 1\right] \\
&= \Psi(s+1) - \Psi\!\left(s + \tfrac{1}{2}\right),
\end{aligned}$$

where the last two steps follow by equation **1.512**(5) in [10] and equation **3.13**(42) in [12], respectively. The proof is complete.